\documentclass[12pt]{article}
\usepackage{graphicx}
\usepackage{amsmath,amsthm,amssymb,enumerate}
\usepackage{euscript,mathrsfs}
\usepackage{color}
\usepackage{dsfont}
\usepackage{url}
\usepackage[left=2cm,right=2cm,top=3.5cm,bottom=3.5cm]{geometry}
\usepackage{color}
\usepackage[framemethod=tikz]{mdframed}
\allowdisplaybreaks

\usepackage{soul}

\catcode`\@=11 \@addtoreset{equation}{section}

\catcode`\@=12

\newtheorem{Theorem}{Theorem}[section]
\newtheorem{Proposition}[Theorem]{Proposition}
\newtheorem{Lemma}[Theorem]{Lemma}
\newtheorem{Corollary}[Theorem]{Corollary}

\theoremstyle{definition}
\newtheorem{Definition}[Theorem]{Definition}

\newtheorem{Remark}[Theorem]{Remark}

\newcommand{\bTheorem}[1]{
	\begin{Theorem} \label{T#1} }
	\newcommand{\eT}{\end{Theorem}}

\newcommand{\bProposition}[1]{
	\begin{Proposition} \label{P#1}}
	\newcommand{\eP}{\end{Proposition}}

\newcommand{\bLemma}[1]{
	\begin{Lemma} \label{L#1} }
	\newcommand{\eL}{\end{Lemma}}

\newcommand{\bCorollary}[1]{
	\begin{Corollary} \label{C#1} }
	\newcommand{\eC}{\end{Corollary}}

\newcommand{\bRemark}[1]{
	\begin{Remark} \label{R#1} }
	\newcommand{\eR}{\end{Remark}}

\newcommand{\bDefinition}[1]{
	\begin{Definition} \label{D#1} }
	\newcommand{\eD}{\end{Definition}}

\newcommand{\Del}{\Delta_x}

\newcommand{\Ds}{\mathbb{D}_x}

\newcommand{\tvm}{\widetilde{\vc{m}}}
\newcommand{\tS}{\widetilde{S}}
\newcommand{\bfphi}{\boldsymbol{\varphi}}

\newcommand{\bFormula}[1]{
	\begin{equation} \label{#1}}
	\newcommand{\eF}{\end{equation}}

\newcommand{\Ov}[1]{\overline{#1}}

\newcommand{\aleq}{\stackrel{<}{\sim}}

\newcommand{\vr}{\varrho}
\newcommand{\vre}{\vr_\ep}

\newcommand{\vte}{\vt_\ep}
\newcommand{\vue}{\vu_\ep}
\newcommand{\tvr}{\tilde \vr}
\newcommand{\tvu}{{\tilde \vu}}
\newcommand{\tvt}{\tilde \vt}

\newcommand{\vt}{\vartheta}
\newcommand{\vu}{\vc{u}}
\newcommand{\vm}{\vc{m}}

\newcommand{\vc}[1]{{\bf #1}}

\newcommand{\Div}{{\rm div}_x}
\newcommand{\Grad}{\nabla_x}

\newcommand{\dx}{\,{\rm d} {x}}

\newcommand{\dt}{\,{\rm d} t }

\newcommand{\vU}{\vc{U}}

\newcommand{\intO}[1]{\int_{\Omega} #1 \ \dx}

\newcommand{\D}{{\rm d}}

\newcommand{\ep}{\varepsilon}

\newcommand{\vtB}{\vt_B}

\newcommand{\br}{ \nonumber \\ }

\def\softd{{\leavevmode\setbox1=\hbox{d}%
		\hbox to 1.05\wd1{d\kern-0.4ex{\char039}\hss}}}
\definecolor{Cgrey}{rgb}{0.85,0.85,0.85}
\definecolor{Cblue}{rgb}{0.50,0.85,0.85}
\definecolor{Cred}{rgb}{1,0,0}
\definecolor{fancy}{rgb}{0.10,0.85,0.10}
\definecolor{amaranth}{rgb}{0.9, 0.17, 0.31}

\newcommand\Cbox[2]{%
	\newbox\contentbox%
	\newbox\bkgdbox%
	\setbox\contentbox\hbox to \hsize{%
		\vtop{
			\kern\columnsep
			\hbox to \hsize{%
				\kern\columnsep%
				\advance\hsize by -2\columnsep%
				\setlength{\textwidth}{\hsize}%
				\vbox{
					\parskip=\baselineskip
					\parindent=0bp
					#2
				}%
				\kern\columnsep%
			}%
			\kern\columnsep%
		}%
	}%
	\setbox\bkgdbox\vbox{
		\color{#1}
		\hrule width  \wd\contentbox %
		height \ht\contentbox %
		depth  \dp\contentbox
		\color{black}
	}%
	\wd\bkgdbox=0bp%
	\vbox{\hbox to \hsize{\box\bkgdbox\box\contentbox}}%
	\vskip\baselineskip%
}

\mdfdefinestyle{MyFrame}{%
	linecolor=black,
	outerlinewidth=1pt,
	roundcorner=5pt,
	innertopmargin=\baselineskip,
	innerbottommargin=\baselineskip,
	innerrightmargin=10pt,
	innerleftmargin=10pt,
	backgroundcolor=white!20!white}

\begin{document}


\title{On the incompressible limit of a strongly stratified 
heat conducting fluid}

\author{Danica Basari\' c $^{1,}$\thanks{The work of D.B. was supported by the
		Czech Sciences Foundation (GA\v CR), Grant Agreement
		21--02411S. The Institute of Mathematics of the Academy of Sciences of
		the Czech Republic is supported by RVO:67985840. } \and Peter Bella $^{2,}$\thanks{P.B. and F.O. were partially supported by Deutsche Forschungsgemeinschaft (DFG) in context of the Emmy Noether Junior Research Group BE 5922/1-1.}
	\and Eduard Feireisl $^{1,}$\thanks{The work of E.F. was partially supported by the
		Czech Sciences Foundation (GA\v CR), Grant Agreement
		21--02411S. The Institute of Mathematics of the Academy of Sciences of
		the Czech Republic is supported by RVO:67985840. } \and Florian Oschmann $^{1,\dagger}$ 
	\and Edriss S. Titi $^{3,}$\thanks{This publication was made possible by NPRP grant \# S-0207-200290 from the Qatar National Research Fund (a member of Qatar Foundation). This research has benefited from the inspiring environment of the CRC 1114 "Scaling Cascades in Complex Systems", Project Number 235221301, Project C06, funded by Deutsche Forschungsgemeinschaft (DFG).}}

\date{}

\maketitle

\medskip

\centerline{$^1$ Institute of Mathematics of the Academy of Sciences of the Czech Republic}

\centerline{\v Zitn\' a 25, CZ-115 67 Praha 1, Czech Republic}

\medskip

\centerline{$^2$ TU Dortmund, Fakultät für Mathematik}

\centerline{Vogelpothsweg 87, 44227 Dortmund, Germany}

\medskip

\centerline{$^3$ Department of Mathematics,  Texas A\&M University, College Station, TX 77843, USA}
	 \centerline{and}
\centerline{Department of Applied Mathematics and Theoretical Physics}
\centerline{University of Cambridge, Cambridge CB3 0WA, UK}

\begin{abstract}
	
	A compressible, viscous and heat conducting fluid is confined between two parallel plates 
	maintained at a constant temperature and subject to a strong stratification due to the gravitational force. 
	We consider the asymptotic limit, where the Mach number and the Froude number are of the same order 
	proportional to a small parameter. We show the limit problem can be identified with Majda's model of 
	layered ``stack-of-pancake'' flow.

\end{abstract}




\section{Introduction}
\label{i}

Consider the motion of a compressible viscous and heat conducting fluid
confined between two parallel plates. For simplicity, we suppose the motion is space--periodic with respect to the horizontal variable. Consequently, the spatial domain $\Omega$ may be identified with 
\[
\Omega = \mathbb{T}^{2} \times (0,1),\ \mathbb{T}^{2} = \left( [-1,1] \Big|_{\{ -1,1 \} } \right)^{2}.
\]
The time evolution of the fluid 
mass density $\vr = \vr(t,x)$, the absolute temperature $\vt = \vt(t,x)$, and the velocity $\vu = \vu(t,x)$ 
is governed by the \emph{Navier--Stokes--Fourier (NSF) system}: 
\begin{align} 
	\partial_t \vr + \Div (\vr \vu) &= 0, \label{i1}\\
	\partial_t (\vr \vu) + \Div (\vr \vu \otimes \vu) + \frac{1}{\ep^2} \Grad p(\vr, \vt) &= \Div \mathbb{S}(\vt, \Grad \vu) + \frac{1}{\ep^{2}} \vr \Grad G, \label{i2} \\ 
	\partial_t (\vr s(\vr, \vt)) + \Div (\vr s (\vr, \vt) \vu) + \Div \left( \frac{ \vc{q} (\vt, \Grad \vt) }{\vt} \right) &= 
	\frac{1}{\vt} \left( \ep^2 \mathbb{S} : \Grad \vu - \frac{\vc{q} (\vt, \Grad \vt) \cdot \Grad \vt }{\vt} \right),
	\label{i3}	
\end{align}
supplemented with the Dirichlet boundary conditions 
\begin{align} 
	\vu|_{\partial \Omega} &= 0, \label{i4} \\
	\vt|_{\partial \Omega} &= \vtB.
	\label{i5}
	\end{align}
The viscous stress tensor is given by Newton's rheological law 
\begin{equation} \label{i6}
	\mathbb{S}(\vt, \Grad \vu) = \mu(\vt) \left( \Grad \vu + \Grad^t \vu - \frac{2}{3} \Div \vu \mathbb{I} \right) + \lambda(\vt) \Div \vu \mathbb{I}, 
	\end{equation}
and the internal energy flux by Fourier's law
\begin{equation} \label{i7} 
	\vc{q}(\vt, \Grad \vt) = - \kappa (\vt) \Grad \vt.
	\end{equation}
The quantity $s = s(\vr, \vt)$ in \eqref{i3} is the entropy of the system, related to the pressure $p = p(\vr, \vt)$ and the internal energy $e = e(\vr, \vt)$ through Gibbs' equation 
\begin{equation} \label{i9} 
	\vt D s = D e + p D \left( \frac{1}{\vr} \right).
	\end{equation}

The potential $G$ represents the effect of gravitation. The
Mach number ${\rm Ma} = \ep$ and the Froude number ${\rm Fr} = \ep$ are both proportional to a small parameter. 
If $\ep > 0$ is small, the fluid is almost incompressible and strongly stratified, cf.~Klein et al.~\cite{KBSMRMHS}.
Our goal is to identify the limit problem for $\ep \to 0$.

\subsection{Asymptotic limit}

In accordance with the scaling of \eqref{i2}, \eqref{i3},
the zero--th order terms in the asymptotic limit are determined by the stationary (static) problem 
\begin{equation} \label{i10}
	\Grad p(\vr, \vt ) =  \vr \Grad G.
\end{equation} 

Applying ${\bf curl}$ operator to identity \eqref{i10}, we successively deduce 
\[
\Grad \vr \times \Grad G = 0, 
\]
and
\[
\frac{\partial p (\vr, \vt)}{\partial \vr} \Grad \vr + 
\frac{\partial p (\vr, \vt)}{\partial \vt} \Grad \vt = \vr \Grad G \ 
\Rightarrow \ \Grad \vt \times \Grad G = 0,
\]
where we have anticipated that the pressure also depends non-trivially on the temperature $\vt$ and is such that $\frac{\partial p(\vr, \vt)}{\partial \vt} \neq 0$. Thus for the static problem \eqref{i10} to be solvable, both $\Grad \vr$ and $\Grad \vt$ must be parallel to $\Grad G$. This fact imposes certain restrictions on the distribution of the boundary temperature $\vtB$. In particular, the motion in an inclined layer studied 
by Daniels et al.~\cite{DBPB} does not admit any static solution. Accordingly, we focus on the particular case 
\begin{align}
G = - g x_3,\ \vtB = \left\{ \begin{array}{l} \Theta_{\rm up} \ \mbox{if}\ x_3 = 1, \\ 
	\Theta_{\rm bott} \ \mbox{if}\ x_3 = 0, \end{array} \right. \br 
\mbox{where}\ g > 0, \  \Theta_{\rm up} > 0, \	\Theta_{\rm bott} > 0 \ \mbox{are constant.}
\label{i17}
\end{align}
Fixing the temperature profile $\vt_B = \Theta(x_3)$ to comply with the boundary conditions \eqref{i17}, we may recover $\vr = r(x_3)$ as a solution of the ODE 
\begin{equation} \label{i18}
\frac{\partial p (r, \Theta)}{\partial \vr} \partial_{x_3} r + \frac{\partial p (r, \Theta)}{\partial \vt} \partial_{x_3} \Theta = - 
r g.
\end{equation}
Needless to say, such a problem may admit infinitely many solutions.

To simplify, we focus on the case $\Theta_{\rm bott} = \Theta_{\rm up} > 0$. Accordingly, we consider 
the reference temperature profile $\Theta = \Theta_{\rm bott} = \Theta_{\rm up}$ -  a positive constant. Then it follows from \eqref{i18} that the static density profile $r = r(x_3)$ must be non--constant as long as $g \ne 0$. Anticipating the asymptotic limit 
\[
\vre \to r ,\ \vte \to \Theta,\ \vue \to \vU \ \mbox{(in some sense)}
\] 
we deduce from the equation of continuity \eqref{i1} 
\begin{equation} \label{i19}
	\Div (r \vU ) = 0.
\end{equation}	
Applying (formally) the same argument to the entropy balance \eqref{i3} we get 
\begin{equation} \label{i20} 
	\Div (r s(r, \Theta) \vU ) = 0.
	\end{equation}
Equations \eqref{i19}, \eqref{i20} are compatible only if 
\[
\Grad r \cdot \vU = 0. 
\]
As $r$ depends only on the vertical $x_3$-variable, this yields 
\begin{equation} \label{i21} 
	U_3 \equiv 0 .
	\end{equation}

In view of the previous arguments, the limit fluid motion exhibits the ``stack of pancakes structure'' described in Chapter 6 of Majda's book \cite{MAJ3}. Specifically, $\vU = [\vU_h, 0]$, and
\begin{align} 
\frac{\partial p (r, \Theta)}{\partial \vr} \partial_{x_3} r &= - rg \label{i22}, \\ 
{\rm div}_h \vU_h &= 0 \label{i23}, \\
r \Big( \partial_t \vU_h + \vU_h \cdot \nabla_h \vU_h \Big) + \nabla_h \Pi &= \mu(\Theta) \Delta_h \vU_h + \mu(\Theta) \partial^2_{x_3, x_3} \vU_h. 	
	\label{i24}	
	\end{align}
Here and hereafter, the subscript $h$ refers to the 
horizontal variable $x_h = (x_1,x_2)$, $\nabla_h =[\partial_{x_1}, \partial_{x_2}]$, ${\rm div}_h \vc{v} = 
\nabla_h \cdot \vc{v}$, $\Delta_h = {\rm div}_h \nabla_h$. The fluid motion is purely horizontal, the coupling between different layers only through the vertical component of the viscous stress. 
	
To the best of our knowledge, there is no rigorous justification of the system \eqref{i22}--\eqref{i24} available in the literature except the inviscid case discussed in \cite{FeKlKrMa}. It is worth noting that a similar problem 
for the barotropic Navier--Stokes system gives rise to a \emph{different} limit, namely the so--called anelastic approximation, see Masmoudi \cite{MAS2} or Feireisl et al.~\cite{FMNS}. Furthermore, as observed in \cite{BelFeiOsc}, the related case of a \emph{low stratification} with $\rm{Ma}=\ep$ and $\rm{Fr} = \sqrt{\ep}$ leads to a limiting system of Oberbeck-Boussinesq type with non-local boundary conditions for the temperature.

\subsection{The strategy of the convergence proof}

We start with the concept of \emph{weak} solutions for the NSF system with Dirichlet boundary conditions 
introduced in \cite{ChauFei}. In particular, we recall the ballistic energy and the associated relative energy 
inequality in Section \ref{L}. Next, we introduce the concept of \emph{strong} solutions to Majda's system 
in Section \ref{S}. In Section \ref{M}, we state our main result. 

The strategy is of type ``weak'' $\to$ ``strong'', meaning the strong solution of the target system is used as a ``test function'' in the relative energy inequality associated to the primitive system.
In Section \ref{B}, we derive the basic energy estimates that control the amplitude of the fluid velocity as well as 
the distance of the density and temperature profiles from their limit values independent of the scaling 
parameter $\ep$.
In Section \ref{co}, we show convergence to the target system \eqref{i22}--\eqref{i24} anticipating the latter admits a regular solution. This formal argument is made rigorous in Section~\ref{E}, where global existence 
for Majda's model is established. The last result may be of independent interest. 

\section{Weak solutions to the primitive NSF system}
\label{L}

Our analysis is based on the concept of weak solutions to the NSF system
introduced in \cite{ChauFei}, cf.~also \cite{FeiNovOpen}.

\begin{Definition}[{\bf Weak solution to the NSF system}] \label{DL1}
	We say that a trio $(\vr, \vt, \vu)$ is a weak solution of the NSF system \eqref{i1}--\eqref{i7}, 
	with the initial data
	\[
	\vr(0, \cdot) = \vr_0,\ \vr \vu (0, \cdot) = \vr_0 \vu_0,\ 
	\vr s(0, \cdot) = \vr_0 s(\vr_0, \vt_0),
	\]
if the following holds:	

\begin{itemize}
	
	\item The solution belongs to the {\bf regularity class}: 
	\begin{align}
		\vr &\in L^\infty(0,T; L^\gamma(\Omega)) \ \mbox{for some}\ \gamma > 1,\ \vr \geq 0 
		\ \mbox{a.a.~in}\ (0,T) \times \Omega, \br
		\vu &\in L^2(0,T; W^{1,2}_0 (\Omega; R^3)), \br 
		\vt^{\beta/2} ,\ \log(\vt) &\in L^2(0,T; W^{1,2}(\Omega)) \ \mbox{for some}\ \beta \geq 2,\ 
		\vt > 0 \ \mbox{a.a.~in}\ (0,T) \times \Omega, \br
		(\vt - \vtB) &\in L^2(0,T; W^{1,2}_0 (\Omega)),
		\label{Lw6}
	\end{align}
where $\vt_B$ is an extension of the boundary data to the whole $\Omega$.	 
	\item The {\bf equation of continuity} \eqref{i1} is satisfied in the sense of distributions, 
	\begin{align} 
		\int_0^T \intO{ \Big[ \vr \partial_t \varphi + \vr \vu \cdot \Grad \varphi \Big] } \dt &=  - 
		\intO{ \vr(0) \varphi(0, \cdot) }
		\label{Lw4}
	\end{align}
	for any $\varphi \in C^1_c([0,T) \times \Ov{\Omega} )$.
	\item The {\bf momentum equation} \eqref{i2} is satisfied in the sense of distributions, 
	\begin{align}
		\int_0^T &\intO{ \left[ \vr \vu \cdot \partial_t \bfphi + \vr \vu \otimes \vu : \Grad \bfphi + 
			\frac{1}{\ep^2} p(\vr, \vt) \Div \bfphi \right] } \dt \br &= \int_0^T \intO{ \left[ \mathbb{S}(\vt, \Grad \vu) : \Grad \bfphi - \frac{1}{\ep^{2}} \vr \Grad G \cdot \bfphi \right] } \dt - 
		\intO{ \vr_0 \vu_0 \cdot \bfphi (0, \cdot) }
		\label{Lw5}
	\end{align}	
	for any $\bfphi \in C^1_c([0, T) \times \Omega; R^3)$.
	
	\item The {\bf entropy balance} \eqref{i3} is replaced by the inequality
	\begin{align}
		- \int_0^T &\intO{ \left[ \vr s(\vr, \vt) \partial_t \varphi + \vr s (\vr ,\vt) \vu \cdot \Grad \varphi + \frac{\vc{q} (\vt, \Grad \vt )}{\vt} \cdot 
			\Grad \varphi \right] } \dt \br &\geq \int_0^T \intO{ \frac{\varphi}{\vt} \left[ \ep^2 \mathbb{S}(\vt, \Grad \vu) : \Ds \vu - 
			\frac{\vc{q} (\vt, \Grad \vt) \cdot \Grad \vt }{\vt} \right] } \dt + \intO{ \vr_0 s(\vr_0, \vt_0) 
			\varphi (0, \cdot) } 
		\label{Lw7} 
	\end{align}
	for any $\varphi \in C^1_c([0, T) \times \Omega)$, $\varphi \geq 0$, where $\Ds \vu = \frac12 ( \Grad \vu + \Grad^t \vu)$ is the symmetric gradient.
	
	\item  The {\bf ballistic energy balance}
	\begin{align}  
		- \int_0^T &\partial_t \psi	\intO{ \left[ \ep^2 \frac{1}{2} \vr |\vu|^2 + \vr e(\vr, \vt) - \tvt \vr s(\vr, \vt) \right] } \dt  \br &+ \int_0^T \psi
		\intO{ \frac{\tvt}{\vt}	 \left[ \ep^2 \mathbb{S}(\vt, \Grad \vu): \Ds \vu - \frac{\vc{q}(\vt, \Grad \vt) \cdot \Grad \vt }{\vt} \right] } \dt  \br
		&\leq 
		\int_0^T \psi \intO{ \left[ \vr \vu \cdot \Grad G - 
			\vr s(\vr, \vt) \partial_t \tvt - \vr s(\vr, \vt) \vu \cdot \Grad \tvt - \frac{\vc{q}(\vt, \Grad \vt)}{\vt} \cdot \Grad \tvt \right] } \dt \br 
		&+ \psi(0) \intO{ \left[ \frac{1}{2} \ep^2 \vr_0 |\vu_0|^2 + \vr_0 e(\vr_0, \vt_0) - \tvt(0, \cdot) \vr_0 s(\vr_0, \vt_0) \right] }
		\label{Lw8}
	\end{align}
	holds for any $\psi \in C^1_c ([0, T))$, $\psi \geq 0$, and any $\tvt \in C^1([0, T) \times \Ov{\Omega})$, 
	\[
	\tvt > 0,\ \tvt|_{\partial \Omega} = \vtB.
	\]
\end{itemize}
 
	\end{Definition}

\subsection{Relative energy inequality}

In addition to Gibbs' equation \eqref{i9}, we impose the hypothesis of thermodynamic stability written in the form 
\begin{equation} \label{HTS}
	\frac{\partial p(\vr, \vt) }{\partial \vr } > 0,\ 
	\frac{\partial e(\vr, \vt) }{\partial \vt } > 0 \ \mbox{for all}\ \vr, \vt > 0.
\end{equation}

Next, following \cite{ChauFei}, we introduce the scaled \emph{relative energy} 
\begin{align}
	E_\ep &\left( \vr, \vt, \vu \Big| \tvr , \tvt, \tvu \right) \br &= \frac{1}{2}\vr |\vu - \tvu|^2 + 
	\frac{1}{\ep^2} \left[ \vr e - \tvt \Big(\vr s - \tvr s(\tvr, \tvt) \Big)- 
	\Big( e(\tvr, \tvt) - \tvt s(\tvr, \tvt) + \frac{p(\tvr, \tvt)}{\tvr} \Big)
	(\vr - \tvr) - \tvr e (\tvr, \tvt) \right] .
	\nonumber
\end{align}
Now, the hypothesis of thermodynamic stability \eqref{HTS} can be equivalently rephrased as (strict) convexity 
of the total energy expressed with respect to the conservative entropy variables
\[
E_\ep \Big( \vr, S = \vr s(\vr, \vt), \vm = \vr \vu \Big) \equiv \frac{1}{2} \frac{|\vm|^2}{\vr} + 
\frac{1}{\ep^2} \vr e(\vr, S),
\]
whereas the relative energy can be written as
\begin{align}
E_\ep &\left( \vr, S, \vm \Big| \tvr , \tS, \tvm \right) = E_\ep(\vr, S, \vm) - \left< \partial_{\vr, S, \vm} E_\ep(\tvr, \tS, \tvm) ; (\vr - \tvr, S - \tS, \vm - \tvm) \right> - E_\ep(\tvr, \tS, \tvm).
\nonumber
\end{align}

Finally, as observed in \cite{ChauFei}, any weak solution in the sense of Definition \ref{DL1} satisfies the \emph{relative energy inequality} 
 \begin{align}
	&\left[ \intO{ E_\ep \left(\vr, \vt, \vu \Big| \tvr, \tvt, \tvu \right) } \right]_{t = 0}^{t = \tau} \br 
	&+ \int_0^\tau \intO{ \frac{\tvt}{\vt} \left( \mathbb{S} (\vt, \Grad \vu) : \Ds \vu + \frac{1}{\ep^2} \frac{\kappa(\vt) |\Grad \vt|^2 }{\vt} \right) } \dt \br 
	&\leq - \frac{1}{\ep^2} \int_0^\tau \intO{ \left( \vr (s - s(\tvr, \tvt)) \partial_t \tvt + \vr (s - s(\tvr, \tvt)) \vu \cdot \Grad \tvt -
		\left( \frac{\kappa (\vt) \Grad \vt}{\vt} \right) \cdot \Grad \tvt \right) } \dt \br 
	&- \int_0^\tau \intO{ \Big[ \vr (\vu - \tvu) \otimes (\vu - \tvu) + \frac{1}{\ep^2} p(\vr, \vt) \mathbb{I} - \mathbb{S}(\vt, \Grad \vu) \Big] : \Ds \tvu } \dt \br 
	&+ \int_0^\tau \intO{ \vr \left[ \frac{1}{\ep^{2}} \Grad G  - \partial_t \tvu - (\tvu \cdot \Grad) \tvu  \right] \cdot (\vu - \tvu) } \dt \br 
	&+ \frac{1}{\ep^2} \int_0^\tau \intO{ \left[ \left( 1 - \frac{\vr}{\tvr} \right) \partial_t p(\tvr, \tvt) - \frac{\vr}{\tvr} \vu \cdot \Grad p(\tvr, \tvt) \right] } \dt
	\label{L4}
\end{align}
for a.a. $\tau > 0$ and any trio of continuously differentiable functions $(\tvr, \tvt, \tvu)$ satisfying
\begin{equation} \label{L5}
	\tvr > 0,\ \tvt > 0,\ \tvt|_{\partial \Omega} = \vtB, \ \tvu|_{\partial \Omega} = 0.
\end{equation}

\subsection{Constitutive relations}
\label{CR}

The existence theory developed in \cite{ChauFei} is conditioned by certain restrictions imposed on the 
constitutive relations (state equations) similar to those introduced in the monograph
\cite[Chapters 1,2]{FeNo6A}. Specifically, the equation of state reads
\[
p(\vr, \vt) = p_{\rm m} (\vr, \vt) + p_{\rm rad}(\vt), 
\]
where $p_{\rm m}$ is the pressure of a general \emph{monoatomic} gas, 
\begin{equation} \label{con1}
	p_{\rm m} (\vr, \vt) = \frac{2}{3} \vr e_{\rm m}(\vr, \vt),
\end{equation}
enhanced by the radiation pressure 
\[
p_{\rm rad}(\vt) = \frac{a}{3} \vt^4,\ a > 0.
\]
Accordingly, the internal energy reads 
\[
e(\vr, \vt) = e_{\rm m}(\vr, \vt) + e_{\rm rad}(\vr, \vt),\ e_{\rm rad}(\vr, \vt) = \frac{a}{\vr} \vt^4.
\]
Moreover, using several physical principles it was shown in \cite[Chapter 1]{FeNo6A}:

\begin{itemize}
	
	\item {\bf Gibbs' relation} together with \eqref{con1} yield 
	\[
	p_{\rm m} (\vr, \vt) = \vt^{\frac{5}{2}} P \left( \frac{\vr}{\vt^{\frac{3}{2}}  } \right)
	\]
	for a certain $P \in C^1[0,\infty)$.
	Consequently, 
	\begin{equation} \label{w9}
		p(\vr, \vt) = \vt^{\frac{5}{2}} P \left( \frac{\vr}{\vt^{\frac{3}{2}}  } \right) + \frac{a}{3} \vt^4,\ 
		e(\vr, \vt) = \frac{3}{2} \frac{\vt^{\frac{5}{2}} }{\vr} P \left( \frac{\vr}{\vt^{\frac{3}{2}}  } \right) + \frac{a}{\vr} \vt^4, \ a > 0.
	\end{equation}
	
	\item {\bf Hypothesis of thermodynamic stability} \eqref{HTS} expressed in terms of 
	$P$ gives rise to
	\begin{equation} \label{w10}
		P(0) = 0,\ P'(Z) > 0 \ \mbox{for}\ Z \geq 0,\ 0 < \frac{ \frac{5}{3} P(Z) - P'(Z) Z }{Z} \leq c \ \mbox{for}\ Z > 0.
	\end{equation} 	
	In particular, the function $Z \mapsto P(Z)/ Z^{\frac{5}{3}}$ is decreasing, and we suppose 
	\begin{equation} \label{w11}
		\lim_{Z \to \infty} \frac{ P(Z) }{Z^{\frac{5}{3}}} = p_\infty > 0.
	\end{equation}
	
	\item 
	Accordingly, the associated entropy
	takes the form 
	\begin{equation} \label{w12}
		s(\vr, \vt) = s_{\rm m}(\vr, \vt) + s_{\rm rad}(\vr, \vt),\ s_{\rm m} (\vr, \vt) = \mathcal{S} \left( \frac{\vr}{\vt^{\frac{3}{2}} } \right),\ s_{\rm rad}(\vr, \vt) = \frac{4a}{3} \frac{\vt^3}{\vr}, 
	\end{equation}
	where 
	\begin{equation} \label{w13}
		\mathcal{S}'(Z) = -\frac{3}{2} \frac{ \frac{5}{3} P(Z) - P'(Z) Z }{Z^2} < 0.
	\end{equation}
	In addition, we impose the {\bf Third law of thermodynamics}, cf.~Belgiorno \cite{BEL1}, \cite{BEL2}, requiring the entropy to vanish 
	when the absolute temperature approaches zero, 
	\begin{equation} \label{w14}
		\lim_{Z \to \infty} \mathcal{S}(Z) = 0.
	\end{equation}
	
\end{itemize}

Finally, we
suppose the transport coefficients are continuously differentiable functions satisfying
\begin{align} 
	0 < \underline{\mu}(1 + \vt) &\leq \mu(\vt),\ |\mu'(\vt)| \leq \Ov{\mu}, \br 
	0 &\leq \eta (\vt) \leq \Ov{\eta}(1 + \vt), \br
	0 < \underline{\kappa} (1 + \vt^\beta) &\leq \kappa (\vt) \leq \Ov{\kappa}(1 + \vt^\beta), 
	\ \mbox{where}\ \beta > 6. \label{w16}
\end{align}

As a consequence of the above hypotheses, we get the following estimates:
\begin{align} 
	\vr^{\frac{5}{3}} + \vt^4 \aleq \vr e(\vr, \vt) &\aleq 	1+ \vr^{\frac{5}{3}} + \vt^4, \label{L5b} \\
	s_{\rm m}(\vr, \vt) &\aleq \left( 1 + |\log(\vr)| + [\log(\vt)]^+ \right), \label{L5a}
\end{align} 
see \cite[Chapter 3, Section 3.2]{FeNo6A}.

\section{Strong solutions to Majda's system}
\label{S}

Problem \eqref{i23}--\eqref{i24} shares many common features with the $2d-$incompressible Navier--Stokes system solved in the celebrated work by Lady\v zenskaja \cite{LAD5}, \cite{LAD6}. Indeed we show that problem 
\eqref{i23}--\eqref{i24}, endowed with the boundary conditions 
\begin{equation} \label{S1}
\vU_h |_{\partial \Omega} = 0,\ \Omega = \mathbb{T}^2 \times (0,1),\ 
\mathbb{T}^2 = \left( [-1,1] \Big|_{\{ -1,1 \} } \right)^{2},	
	\end{equation}
is globally well posed in the framework of Sobolev spaces $W^{2,p}$ with $p > 1$ large enough. We report the following result that may be of independent interest.

\begin{mdframed}[style=MyFrame]  

	\begin{Theorem}[{\bf Global existence for Majda's system}] \label{TH1}
	
	Let $\Theta > 0$ be given. Suppose that 
	\begin{equation} \label{H1}
		r \in C^1([0,1]),\ 0 < \underline{r} \leq r (x_3) \ \mbox{for all}\ x_3 \in [0,1]. 
	\end{equation}
	Let the initial data $\vU_{0,h}$ belong to the class 
	\begin{equation} \label{H2}
		\vU_{0,h} \in W^{3,q} \cap W^{1,q}_0(\Omega; R^2),\ {\rm div}_h \vU_{0,h} = 0 
	\end{equation}
	for all $1 \leq q < \infty$. 
	
	Then the system \eqref{i23}--\eqref{i24}, with the boundary conditions \eqref{S1} and the initial condition \eqref{H2}, admits a strong solution 
	$\vU_h$ in $(0,T) \times \Omega$, unique in the class 
	\begin{equation} \label{H3}
		\partial_t \vU_h \in L^p(0,T; L^p(\Omega; R^2)),\ 
		(\vU_h , \nabla_h \vU_h) \in L^p(0,T; W^{2,p}(\Omega; R^2)\times  W^{2,p}(\Omega; R^{2\times 2} ))	
	\end{equation}		
	for any $1 \leq p < \infty$.	 
\end{Theorem}

\end{mdframed}

\begin{Remark}
	To avoid any misunderstanding we emphasize that by 
\[
	\vU_{0,h} \in W^{3,q} \cap W^{1,q}_0(\Omega; R^2),\ {\rm div}_h \vU_{0,h} = 0 
\]
for all $1 \leq q < \infty$ we mean
\[
\vU_{0,h} \in \bigcap_{q \geq 1} W^{3,q} \cap W^{1,q}_0(\Omega; R^2),\ {\rm div}_h \vU_{0,h} = 0. 
\]
Similarly, 
\[
	\partial_t \vU_h \in L^p(0,T; L^p(\Omega; R^2)),\ 
	(\vU_h , \nabla_h \vU_h) \in L^p(0,T; W^{2,p}(\Omega; R^2)\times  W^{2,p}(\Omega; R^{2\times 2} ))	
\]		
for all finite $1 \leq p < \infty$ means
\[
\partial_t \vU_h \in \bigcap_{p \geq 1} L^p(0,T; L^p(\Omega; R^2)),\ 
(\vU_h , \nabla_h \vU_h) \in \bigcap_{p \geq 1} L^p(0,T; W^{2,p}(\Omega; R^2)\times  W^{2,p}(\Omega; R^{2\times 2} )).	
\]

	\end{Remark}

The proof of Theorem \ref{TH1} is postponed to Section \ref{E}.

\section{Main result}
\label{M}

Having collected the necessary preliminary material, we are ready to state our main result. 

\begin{mdframed}[style=MyFrame]

\begin{Theorem}[{\bf Singular limit}] \label{MT1}
	
	Let the thermodynamic functions $p$, $e$, and $s$ as well as the transport coefficients $\mu$, $\lambda$, and $\kappa$ comply with the structural hypotheses specified in Section~\ref{CR}. Let 
	\begin{equation} \label{S2}
	G = -g x_3,\ g > 0,\ \Theta_{\rm up} = \Theta_{\rm bott} = \Theta > 0, 	
		\end{equation}
and let
\begin{equation} \label{S3}
r \in C^1([0,1]) ,\ 0 < \underline{r} \leq r,\ \frac{\partial p(r, \Theta) }{\partial \vr} \partial_{x_3} r = - r g.	
	\end{equation}
Let $(\vre, \vte, \vue)_{\ep > 0}$ be a family of weak solutions of the scaled NSF system in the sense of Definition \ref{DL1} emanating from the initial data 
\[
\vre(0,\cdot) = \vr_{0,\ep},\ \vre \vue (0,\cdot) = \vr_{0,\ep} \vu_{0,\ep},\ \vre s(\vre, \vte)(0, \cdot) = \vr_{0,\ep} s(\vr_{0,\ep}, \vt_{0,\ep}), 
\]
where 
\begin{equation} \label{S4}
	\intO{ E_\ep \left( \vr_{0,\ep}, \vt_{0,\ep}, \vu_{0, \ep}\ \Big| \ r, \Theta, [\vU_{0,h} , 0 ] \right) } \to 0 \ \mbox{as}\ \ep \to 0, 
\end{equation}
and $\vU_{0,h}$ belongs to the class \eqref{H2}. 

Then 
\begin{equation} \label{S5}
{\rm ess} \sup_{\tau \in (0,T)}	
\intO{ E_\ep \left( \vr_{\ep}, \vt_{\ep}, \vu_{\ep}\ \Big| \ r, \Theta, [ \vU_{h}, 0 ] \right) (\tau, \cdot) } \to 0 \ \mbox{as}\ \ep \to 0,	
	\end{equation}
where $\vU_h$ is the unique solution of Majda's system, the existence of which is guaranteed by Theorem \ref{TH1}.

	\end{Theorem}

\end{mdframed}

Hypothesis \eqref{S4} corresponds to \emph{well--prepared} initial data. In view of the coercivity properties of the relative energy stated in \eqref{BB1}, \eqref{BB2} below, relation 
\eqref{S5} implies, in particular, 
\begin{align*}
	\vre &\to r \hspace{1.4cm} \mbox{in}\ L^\infty(0,T; L^{\frac{5}{3}}(\Omega)), \\
	\vte &\to \Theta \hspace{1.3cm} \mbox{in}\ L^\infty(0,T; L^2(\Omega)), \\
	\vre \vue &\to r [\vU_h, 0] \ \ \mbox{in}\ L^\infty(0,T; L^{1}(\Omega; R^3))
\end{align*}
as $\ep \to 0$.

The next two sections are devoted to the proof of Theorem \ref{MT1}.

\section{Uniform bounds}
\label{B}

In order to perform the singular limit in the NSF system we need the associated sequence of weak solutions 
$(\vre, \vte, \vue)_{\ep > 0}$ to be bounded at least in the energy space. First, we  
introduce the notation of \cite{FeNo6A} to distinguish between the ``essential'' and ``residual'' range of the 
thermostatic variables $(\vr, \vt)$. Specifically, given a compact set 
\[
K \subset \left\{ (\vr, \vt) \in R^2 \ \Big| \ \vr > 0, \vt > 0 \right\}
\]
we introduce 
\[
g_{\rm ess} = g \mathds{1}_{(\vr, \vt) \in K},\ 
g_{\rm res} = g - g_{\rm ess} = g \mathds{1}_{(\vr, \vt) \in R^2 \setminus K}.
\]
As shown in \cite[Chapter 5, Lemma 5.1]{FeNo6A}, the relative energy enjoys the following coercivity properties:
\begin{align} 
E_{\ep} \left( \vr, \vt, \vu \Big| \tvr, \tvt, \tvu \right)	&\geq E_{\ep} \left( \vr, \vt, \vu \Big| \tvr, \tvt, \tvu \right)_{\rm ess} \geq 
C \left( \frac{ |\vr - \tvr|^2 }{\ep^2} + \frac{ |\vt - \tvt|^2 }{\ep^2} + |\vu - \tvu |^2 \right)_{\rm ess} 
	\label{BB1} \\ 
E_{\ep} \left( \vr, \vt, \vu \Big| \tvr, \tvt, \tvu \right)	&\geq E_{\ep} \left( \vr, \vt, \vu \Big| \tvr, \tvt, \tvu \right)_{\rm res} \geq 
C  \left( \frac{1}{\ep^2} + \frac{1}{\ep^2} \vr e(\vr, \vt) + \frac{1}{\ep^2} \vr |s(\vr, \vt)| + \vr |\vu|^2 \right)_{\rm res} 
\label{BB2}	
	\end{align}
whenever $(\tvr, \tvt) \in {\rm int}[K]$, where the constant $C$ depends on $K$ and the distance 
\[
{\rm dist} \left[ (\tvr , \tvt  ) ; \partial K \right]. 
\]

\subsection{Energy estimates for ill--prepared data}
\label{ee}

We examine a slightly more general situation than in Theorem \ref{MT1}. Let $\Theta>0$ be constant and $r$ the solution of the static problem 
\begin{equation} \label{SP}
	\frac{ \partial p(r, \Theta )}{\partial \vr} \partial_{x_3} r = - r g. 
\end{equation} 

Next, we consider a family $(\vre, \vte, \vue)_{\ep > 0}$ emanating from \emph{ill--prepared} data 
$(\vr_{0,\ep}, \vt_{0,\ep}, \vu_{0, \ep})_{\ep > 0}$,
\begin{equation} \label{SL1}
	\intO{ E_\ep \left( \vr_{0, \ep}, \vt_{0,\ep}, \vu_{0, \ep} \Big| r , \Theta , 0 \right) } \aleq 1 \ \mbox{independently of}\ \ep \to 0.
\end{equation}
The relative energy inequality \eqref{L4} yields 
\begin{align}
	&\left[ \intO{ E_\ep \left(\vre, \vte, \vue \Big| r, \Theta , 0 \right) } \right]_{t = 0}^{t = \tau} \br 
	&+ \int_0^\tau \intO{ \frac{\Theta }{\vte} \left( \mathbb{S} (\vte, \Grad \vue) : \Ds \vue + \frac{1}{\ep^2} \frac{\kappa (\vte) |\Grad \vte |^2 }{\vte} \right) } \dt \br 
	&\leq \frac{1}{\ep^2} \int_0^\tau \intO{ \frac{\vre}{r} \left( r\Grad G -  \Grad p(r, \Theta ) \right)  \cdot \vue  } \dt.
	\label{SL2}
\end{align}
Moreover, in view of \eqref{S2} and \eqref{S3}, we deduce the stationary equation
\begin{equation}\label{stationary equation}
  \Grad p(r, \Theta ) = r \Grad G;
\end{equation}
hence \eqref{SL2} reduces to 
\begin{align}
	&\left[ \intO{ E_\ep \left(\vre, \vte, \vue \Big| r, \Theta , 0 \right) } \right]_{t = 0}^{t = \tau} \br 
	&\quad + \int_0^\tau \intO{ \frac{\Theta }{\vte} \left( \mathbb{S} (\vte, \Grad \vue) : \Ds \vue + \frac{1}{\ep^2} \frac{\kappa (\vte) |\Grad \vte |^2 }{\vte} \right) } \dt \leq 0.
	\label{SL2a}
\end{align}

\subsection{Conclusion, uniform bounds for ill-prepared data}
\label{cbe}

In view of the estimates obtained in the previous section, we deduce from \eqref{SL2a} for ill--prepared initial data satisfying \eqref{SL1} the following bounds independent of the scaling parameter $\ep \to 0$: 

\begin{align}
	{\rm ess} \sup_{t \in (0,T)} \intO{ E_\ep \left( \vre, \vte, \vue \Big| r, \Theta , 0 \right) } &\aleq 1, \label{be1} \\
	\int_0^T \| \vue \|^2_{W^{1,2}_0 (\Omega; R^3) } \dt &\aleq 1, \label{be2} \\
	\frac{1}{\ep^2} \int_0^T \left( \| \Grad \log (\vte) \|^2_{L^2(\Omega; R^3)} + \| \Grad \vte^{\frac{\beta}{2}} \|^2_{L^2(\Omega; R^3)} \right) &\aleq 1.
	\label{be3}
	\end{align}

Next, it follows from \eqref{be1} that the measure of the residual set shrinks to zero, specifically
\begin{equation} \label{be4}
	\frac{1}{\ep^2} {\rm ess} \sup_{t \in (0,T)} \intO{ [1]_{\rm res} } \aleq 1.
\end{equation}

In addition, we get from \eqref{be1}:

\begin{align}
	{\rm ess} \sup_{t \in (0,T)} \intO{ \vre |\vue|^2 } &\aleq 1, \br 
	{\rm ess} \sup_{t \in (0,T)} \left\| \left[ \frac{\vre - r}{\ep} \right]_{\rm ess} \right\|_{L^2(\Omega)} &\aleq 1, \br
	{\rm ess} \sup_{t \in (0,T)} \left\| \left[ \frac{\vte - \Theta}{\ep} \right]_{\rm ess} \right\|_{L^2(\Omega)} &\aleq 1, \br
	\frac{1}{\ep^2} {\rm ess} \sup_{t \in (0,T)} \| [\vre]_{\rm res} \|^{\frac{5}{3}}_{L^{\frac{5}{3}}(\Omega)} + \frac{1}{\ep^2}  {\rm ess} \sup_{t \in (0,T)} \| [\vte]_{\rm res} \|^{4}_{L^{4}(\Omega)}	 &\aleq 1.
	\label{be5}	
\end{align}

Combining \eqref{be3}, \eqref{be4}, and \eqref{be5}, we conclude 
\begin{equation} \label{be6}
	\int_0^T \left\| \frac{\log(\vte) - \log(\Theta)}{\ep} \right\|^2_{W^{1,2}(\Omega)} \dt + \int_0^T \left\| \frac{ \vte - \Theta }{\ep} \right\|^2_{W^{1,2}(\Omega)} \dt \aleq 1.
	\end{equation}

Finally, we claim the bound on the entropy flux
\begin{equation} \label{be7} 
\int_0^T \left\| \left[ \frac{\kappa (\vte) }{\vte} \right]_{\rm res} \frac{\Grad \vte }{\ep} \right\|^q_{L^q(\Omega; R^3)} \dt \aleq 1 \ \mbox{for some}\ q > 1.
\end{equation} 
Indeed we have 
\[
\left|  \left[ \frac{\kappa (\vte) }{\vte} \right]_{\rm res} \frac{\Grad \vte }{\ep} \right| \aleq  \frac{1}{\ep}
\left| \Grad \log (\vte) \right| + \frac{1}{\ep} \left| \left[\vte^{\frac{\beta}{2}} \Grad \vte^{\frac{\beta}{2}}\right]_{\rm res} \right|, 
\]	
where the former term on the right--hand side is controlled via \eqref{be6}. As	for the latter, we deduce from 
\eqref{be3} that  
\[
\left\| \frac{1}{\ep} \Grad \vte^{\frac{\beta}{2}} \right\|_{L^2((0,T) \times \Omega; R^3)} \aleq 1;
\]
hence it is enough to check  
\begin{equation} \label{be8}
\left\| \left[\vte^{\frac{\beta}{2}}\right]_{\rm res} \right\|_{L^r ((0,T) \times \Omega)} \aleq 1 \ \mbox{for some}\ r > 2.
\end{equation}
To see \eqref{be8} first observe that 
\begin{equation} \label{be9}
	{\rm ess} \sup_{t \in (0,T)} \| [\vte]_{\rm res} \|_{L^4(\Omega)} \aleq 1,
\end{equation}
and,  in view of \eqref{be3} and Poincar\' e inequality, 
\[
\left\| 
\vte^{\frac{\beta}{2}} \right\|_{L^2(0,T; L^6(\Omega))} \aleq 1.
\]
Consequently, \eqref{be8} follows by interpolation.

Of course, the above uniform bound remain valid also for the well-prepared initial data considered in Theorem \ref{MT1}.

\section{Convergence to the target system}
\label{co}

We show convergence to the regular solution $\vU_h$ in Majda's system claimed in Theorem \ref{MT1}. To get a lean notation, we will identify the two-dimensional velocity $\vc{U}_h$ with its three-dimensional counterpart $[\vc{U}_h, 0]$.
The ansatz $(\tvr, \tvt, \tvu) = (r, \Theta, \vU_h)$ in the relative energy inequality \eqref{L4} 
yields 
 \begin{align}
	&\left[ \intO{ E_\ep \left(\vre, \vte, \vue \Big| r , \Theta , \vU_h \right) } \right]_{t = 0}^{t = \tau} \br 
	&+ \int_0^\tau \intO{ \frac{\Theta}{\vte} \left( \mathbb{S} (\vte, \Grad \vue) : \Ds \vue + \frac{1}{\ep^2} \frac{\kappa(\vte) |\Grad \vte|^2 }{\vte} \right) } \dt \br 
	&\leq 
	- \int_0^\tau \intO{ \Big[ \vre (\vue - \vU_h) \otimes (\vue - \vU_h) + \frac{1}{\ep^2} p(\vre, \vte) \mathbb{I} - \mathbb{S}(\vte, \Grad \vue) \Big] : \Ds \vU_h } \dt \br 
	&+ \int_0^\tau \intO{ \vre \left[  \partial_t \vU_h + (\vU_h \cdot \Grad) \vU_h  \right] \cdot (\vU_h - \vue) } \dt - \frac{1}{\ep^2} \int_0^\tau \intO{ \vre \Grad G \cdot \vU_h } \dt  , 
	\label{co1}
\end{align}
where we have used the stationary equation
\[
\Grad p (r, \Theta) = r \Grad G.
\] 

Next, seeing that 
\[
\Div \vU_h = 0,\ \Grad G \cdot \vU_h = 0,
\]
we deduce 
 \begin{align}
	&\left[ \intO{ E_\ep \left(\vre, \vte, \vue \Big| r , \Theta , \vU_h \right) } \right]_{t = 0}^{t = \tau} \br 
	&+ \int_0^\tau \intO{ \frac{\Theta}{\vte} \left( \mathbb{S} (\vte, \Grad \vue) : \Ds \vue   + \frac{1}{\ep^2} \frac{\kappa(\vte) |\Grad \vte|^2 }{\vte} \right) } \dt \br 
	&\leq 
	- \int_0^\tau \intO{ \Big[ \vre (\vue - \vU_h) \otimes (\vue - \vU_h)  - \mathbb{S}(\vte, \Grad \vue) \Big] : \Ds \vU_h } \dt \br 
	&+ \int_0^\tau \intO{ \vre \left[  \partial_t \vU_h + (\vU_h \cdot \Grad) \vU_h  \right] \cdot (\vU_h - \vue) } \dt . 
	\label{co2}
\end{align}

Now, in view of the uniform bounds \eqref{be2}, \eqref{be5}, 
\begin{align*}
	\int_0^\tau &\intO{ \vre \left[  \partial_t \vU_h + (\vU_h \cdot \Grad) \vU_h  \right] \cdot (\vU_h - \vue) } \dt \\
	&= \int_0^\tau \intO{ r \left[  \partial_t \vU_h + (\vU_h \cdot \Grad) \vU_h  \right] \cdot (\vU_h - \vue) } \dt
	+ \mathcal{Q}(\ep),
\end{align*}
where $\mathcal{Q}(\ep)$ denotes a generic function with the property $\mathcal{Q}(\ep) \to 0$ as $\ep \to 0$. 

Next, in view of \eqref{be2}, \eqref{be5}, we may assume 
\[
\vre \to r \ \mbox{in}\ L^\infty(0,T; L^{\frac{5}{3}}(\Omega)),\ 
\vue \to \vc{u} \ \mbox{weakly in}\ L^2(0,T; W^{1,2}_0(\Omega)), 
\]
up to a suitable subsequence, where 
\begin{equation} \label{cond1}
\Div (r \vu ) = 0.
\end{equation}
Similarly, using the bounds \eqref{be5}, \eqref{be6} we may perform the limit in the 
entropy inequality \eqref{Lw7} obtaining 
\[
\Div (r s(r, \Theta) \vu) \geq 0.
\]
However, thanks to the no--slip boundary conditions, 
\[
\intO{ \Div (r s(r, \Theta) \vu) } = 0; 
\]
therefore 
\begin{equation} \label{cond2}
\Div (r s(r, \Theta) \vu) = 0.
\end{equation}
Combining \eqref{cond1}, \eqref{cond2} we may infer 
\[
r \frac{ \partial s (r, \Theta) }{\partial \vr} \Grad r \cdot \vu = 0. 
\]
As entropy is given by the constitutive equation \eqref{w12}, \eqref{w13}, 
\[
\frac{ \partial s (r, \Theta) }{\partial \vr} < 0,
\]
and we conclude 
\begin{equation} \label{cond3}
	u_3 = 0,\ {\rm div}_h \vu = 0.
	\end{equation}

Now, 
\begin{align*}
	\int_0^\tau &\intO{ \vre \left[  \partial_t \vU_h + (\vU_h \cdot \Grad) \vU_h  \right] \cdot (\vU_h - \vue) } \dt \\
	&= \int_0^\tau \intO{ r \left[  \partial_t \vU_h + (\vU_h \cdot \Grad) \vU_h  \right] \cdot (\vU_h - \vu) } \dt
	+ \mathcal{Q}(\ep) .
\end{align*}
In addition, 
since $\vU_h$, $\vu$ satisfy \eqref{i24}, \eqref{cond3}, respectively, we obtain 
\begin{align}
\int_0^\tau &\intO{ r \left[  \partial_t \vU_h + (\vU_h \cdot \Grad) \vU_h  \right] \cdot (\vU_h - \vu) } \dt\\
&= 
\int_0^\tau \intO{  \mu(\Theta) \left[  \Delta_h \vU_h + \partial^2_{x_3, x_3} \vU_h  \right] \cdot (\vU_h - \vu) } \dt 
\br &= - \int_0^\tau \intO{ \mathbb{S}(\Theta, \Grad \vU_h ) : \Ds (\vU_h - \vu) } \dt.
\nonumber
\end{align}
Going back to \eqref{co2}, we deduce 
 \begin{align}
	&\left[ \intO{ E_\ep \left(\vre, \vte, \vue \Big| r , \Theta , \vU_h \right) } \right]_{t = 0}^{t = \tau} \br 
	&+ \int_0^\tau \intO{ \frac{\Theta}{\vte} \left( \mathbb{S} (\vte, \Grad \vue) : \Ds \vue   + \frac{1}{\ep^2} \frac{\kappa(\vte) |\Grad \vte|^2 }{\vte} \right) } \dt \br 
	&\leq 
	- \int_0^\tau \intO{ \Big[ \vre (\vue - \vU_h) \otimes (\vue - \vU_h)  - \mathbb{S}(\Theta, \Grad \vu) \Big] : \Ds \vU_h } \dt \br 
	&- \int_0^\tau \intO{ \mathbb{S}(\Theta, \Grad \vU_h ) : \Ds (\vU_h - \vu) } \dt+ \mathcal{Q}(\ep) . 
	\label{co3}
\end{align}

Finally, exploiting weak lower semi--continuity of convex functions, we conclude 
 \begin{align}
	&\left[ \intO{ E_\ep \left(\vre, \vte, \vue \Big| r , \Theta , \vU_h \right) } \right]_{t = 0}^{t = \tau} \br 
	&+ \int_0^\tau \intO{ \left( \mathbb{S} (\Theta, \Grad \vu) - \mathbb{S} (\Theta, \Grad \vU_h \right)  : \left( \Ds \vu - \Ds \vU_h \right) } \dt \br 
	&\leq 
	- \int_0^\tau \intO{ \Big[ \vre (\vue - \vU_h) \otimes (\vue - \vU_h) \Big] : \Ds \vU_h } \dt + \mathcal{Q}(\ep) ,
	\label{co4}
\end{align}
which, applying the standard Gr\"onwall argument, yields the desired convergence as well as $\vu = \vU_h$. 

We have proved Theorem \ref{MT1}.

\section{Global existence for Majda's problem}
\label{E}

Our ultimate goal is to show global existence of strong solutions to Majda's model claimed in Theorem \ref{TH1}.
To this end, it is more convenient to consider the (horizontal) vorticity formulation of \eqref{i23}, \eqref{i24}. With a slight abuse of notation in the definition of $\vc{U}_h$, this formulation reads
\begin{align} 
	\partial_t \omega + \vc{U}_h \cdot \Grad \omega &= \nu \Del \omega, \label{A3} \\ 
	\vc{U}_h &= \left[ \nabla^{\perp}_h \Delta_h^{-1}[ \omega ], 0 \right] \label{A4}, \\ 
	\nu &= \nu(x_3), \label{A5}
\end{align}
with the boundary conditions 
\begin{equation} \label{A6}
	\omega|_{\partial \Omega} = 0, 
\end{equation}
and the initial condition 
\begin{equation} \label{A6i}
	\omega(0, \cdot) = \omega_0. 
	\end{equation} 
Here, $\nu = \frac{\mu(\Theta)}{r}$, and 
\begin{equation} \label{A1}
	\omega  = {\rm curl}_h \vU_h, \ {\rm curl}_h [\vc{v}] = \partial_{x_1} v_2 - \partial_{x_2} v_1.  
\end{equation}
For given $\omega$, the velocity field $\vU_h$ can be recovered via Biot-Savart law: 
\begin{equation} \label{A2}
	\vU_h = \left[ \nabla^{\perp}_h \Delta_h^{-1}[ \omega ], 0 \right],\ \nabla^\perp_h = [ -\partial_{x_2},\partial_{x_1} ].	
\end{equation}

\begin{Remark} \label{RA1} 
	
	Strictly speaking, the velocity $\vU_h$ is determined by \eqref{A2} up to its horizontal average 
	\[
	\Ov{\vU}_h = \int_{\mathbb{T}^2} \vU_h \ \D x_h 
	\] 
	that can be recovered as the unique solution of the parabolic problem 
	\begin{align} 
		r \partial_t \Ov{\vU}_h &= \mu(\Theta) \partial^2_{x_3,x_3} \Ov{\vU}_h \ \mbox{in}\ (0,T) \times (0,1), \br 
		\Ov{\vU}_h |_{x_3 = 0,1} &= 0, \br
		\Ov{\vU}_h(0, \cdot) &= \int_{\mathbb{T}^2} \vU_{0,h} \ \D x_h . \nonumber
		\end{align}

	\end{Remark}

\subsection{Construction via a fixed point argument} 

The desired solution $\omega$ to \eqref{A3}--\eqref{A6i} can be constructed via a simple fixed point argument. Consider the set 
\[
X_M = \left\{ \widetilde{\omega} \in C([0,T] \times \Ov{\Omega}) \ \Big| \ \widetilde \omega|_{\partial \Omega} = 0,\  
\widetilde \omega(0, \cdot) = {\rm curl}_h \vU_{0, h}\ 
,\ \| \widetilde \omega \|_{ C([0,T] \times \Ov{\Omega}) } \leq M
\right\}. 
\]
As the initial velocity $\vU_{0,h}$ belongs to the class \eqref{H2}, the set $X_M$ is a bounded closed convex subset of the Banach space $C([0,T] \times \Ov{\Omega})$. 
Moreover, $X_M$ is non-empty as long as $M$ is large enough to accommodate the initial condition. 

We define a mapping  $\mathcal{T}[\widetilde{\omega}] = \omega$, where $\omega$ is the unique solution of the problem 
\begin{align} 
	\partial_t \omega + b_L ( \widetilde{\vc{U}}_h ) \cdot \Grad \omega &= \nu \Del \omega, \label{C1} \\ 
	\widetilde{\vc{U}}_h &= \left[ \nabla^{\perp}_h \Delta_h^{-1}[ \widetilde{\omega} ], 0 \right] \label{C2}, \\ 
	\omega|_{\partial \Omega} &= 0, \label{C3} \\ 
	\omega(0, \cdot) &= {\rm curl}_h \vU_{0,h} \label{C4} ,
\end{align}
for some cut--off function $b_L$. Specifically, 
\[
b_L(\widetilde{\vU}_h) = [ b_L( \widetilde{U}_h^1), b_L(\widetilde{U}_h^2), 0], 
\]
where 
\[
b_L \in L^\infty(R) \cap C^\infty(R),\ b_L(Z) = Z \ \mbox{whenever}\ |Z| \leq L.
\]

\subsubsection{Maximum principle}

Applying the standard maximum principle, we deduce 
\begin{equation} \label{C6}
\sup_{t \in [0,T]} \| \mathcal{T}[ \widetilde \omega ](t,\cdot)  \|_{C(\Ov{\Omega})} = 
\sup_{t \in [0,T]} \| \omega(t,\cdot)  \|_{C(\Ov{\Omega})} = \| \omega(0, \cdot) \|_{C(\Ov{\Omega})} 
\aleq \| \vU_{0,h} \|_{W^{2,q} (\Omega; R^2)} \ \mbox{as long as}\ q > 3.
\end{equation} 
Note carefully that the bound \eqref{C6} depends solely on the initial data. In particular, it is independent of 
the specific form of the cut--off function $b_L$. 

\subsubsection{Maximal $L^p-L^q$ regularity}
In view of hypothesis \eqref{S3},  
\begin{equation*}
	\nu \in C^1([0,1]),\ 0 < \underline{\nu} \leq \nu(x_3) \ \mbox{for any}\ x_3 \in [0,1].
\end{equation*} 
Consequently, we can apply the maximal $L^p-L^q$ regularity estimates, see, e.g., Denk, Hieber, and Pr\" uss 
\cite{DEHIEPR}, to obtain 
\begin{align} 
	\| \partial_t \omega \|_{L^p(0,T; L^q(\Omega))} + 
	\| \omega \|_{L^p(0,T; W^{2,q}(\Omega))} &\leq 
	c(p,q) \left( \| \omega(0, \cdot) \|_{W^{2,q} \cap W^{1,q}_0(\Omega)} + 
	\| b_L (\widetilde{\vU}_h) \cdot \Grad \omega \|_{L^p(0,T; L^q(\Omega))} \right),\br 
	1 &< p,q < \infty. \label{A11} 
\end{align}
Here 
\[
\| \omega(0, \cdot) \|_{W^{2,q} \cap W^{1,q}_0(\Omega)} \aleq \| \vU_{0,h} \|_{W^{3,q}(\Omega; R^2)}, 
\]
while, by interpolation and \eqref{C6},  
\begin{align}
\| b_L (\widetilde{\vU}_h) \cdot \Grad \omega \|_{L^q(\Omega)} \leq L \|  \Grad \omega \|_{L^q(\Omega; R^3)} 
&\leq L \| \omega \|_{W^{2,q}(\Omega)}^\lambda \| \omega \|_{L^q(\Omega)}^{1 - \lambda} \br
&\leq L c(q) \| \vU_{0,h} \|_{W^{3,q}(\Omega; R^2)}^{1 - \lambda}\| \omega \|_{W^{2,q}(\Omega)}^\lambda, 
\ t \in (0,T)
\nonumber
\end{align}
for some $0 < \lambda < 1$. Consequently, it follows from \eqref{A11} and our hypotheses imposed on the initial data that 
\begin{equation} \label{A12} 
	\| \partial_t \mathcal{T} [ \widetilde{\omega} ] \|_{L^p(0,T; L^q(\Omega))} + 
	\| \mathcal{T}[ \widetilde \omega ] \|_{L^p(0,T; W^{2,q}(\Omega))} \leq 
	c\left(p,q, \| \vU_{0,h} \|_{W^{3,q}(\Omega; R^2)} \right) \left( 1 + L \right)  
\end{equation}
for all finite $p,q$.

\subsection{Fixed point}

It follows from the estimates \eqref{C6}, \eqref{A12} that $\mathcal{T}$ is a compact (continuous) mapping 
of $X_M$ into $X_M$ provided $M$ is large enough, therefore, by means of Tikhonov--Schauder fixed point Theorem, there is a fixed point $\omega \in X_M$ satisfying 
\begin{align} 
	\partial_t \omega + b_L (\vc{U}_h ) \cdot \Grad \omega &= \nu \Del \omega, \br 
	\vc{U}_h &= \left[ \nabla^{\perp}_h \Delta_h^{-1}[ \omega ], 0 \right] \br 
	\omega|_{\partial \Omega} &= 0, \br 
	\omega(0, \cdot) &= {\rm curl}_h \vU_{0,h}. \nonumber
\end{align}

Finally, as $\vU_h$ is given by the Biot--Savart law, we get 
\[
	\sup_{x_3 \in (0,1)} \| \nabla_h \vU_h \|_{L^q(\mathbb{T}^2; R^{2 \times 2})} \leq 
	c(q) \| \omega(0, \cdot) \|_{L^\infty (\Omega)} \ \mbox{uniformly for}\ t \in (0,T) 
	\ \mbox{for any}\ 1 < q < \infty,
\]
in particular
\[ 
	\| \vU_h \|_{L^\infty((0,T) \times \Omega; R^2)} \aleq \| \omega(0, \cdot) \|_{L^\infty(\Omega)} 
	\aleq \| \vU_{0,h} \|_{W^{2,q}(\Omega; R^2)} \ \mbox{as soon as}\ q > 3.
\]
Since this bound is independent of $L$, we may choose $L$ large enough 
so that $b_L (\vc{U}_h ) = \vU_h$ to get the desired conclusion 
\begin{align} 
	\partial_t \omega + \vc{U}_h  \cdot \Grad \omega &= \nu \Del \omega, \br 
	\vc{U}_h &= \left[ \nabla^{\perp}_h \Delta_h^{-1}[ \omega ], 0 \right] \br 
	\omega|_{\partial \Omega} &= 0, \br 
	\omega(0, \cdot) &= {\rm curl}_h \vU_{0,h}. \nonumber
\end{align}

Finally, it is easy to check that the solution is unique in the regularity class \eqref{H3}. As a matter of fact, a more general \emph{weak--strong} uniqueness holds 
that could be shown adapting the above arguments based on the relative energy inequality.

We have proved Theorem \ref{TH1}.

\def\cprime{$'$} \def\ocirc#1{\ifmmode\setbox0=\hbox{$#1$}\dimen0=\ht0
	\advance\dimen0 by1pt\rlap{\hbox to\wd0{\hss\raise\dimen0
			\hbox{\hskip.2em$\scriptscriptstyle\circ$}\hss}}#1\else {\accent"17 #1}\fi}



\end{document}